\documentclass[12pt]{article}
\usepackage{amsthm,amsmath, amssymb} 
\begin{document}
\def\a{\alpha}
\def\b{\beta}
\def\ra{\rightarrow}
\def\hs{\heartsuit}
\def\ss{\spadesuit}
\def\cs{\clubsuit}
\title{Universal Cycles of Classes of Restricted Words}         % Enter your title between curly braces
\author{Arielle Leitner and Anant Godbole}
\date{}          % Enter your date or \today between curly braces
\maketitle

\def\ra{\rightarrow}
\def\ld{\ldots}
\begin{abstract}
It is well known that Universal Cycles (U-Cycles) of $k$-letter words on an $n$-letter alphabet exist for all $k$ and $n$.  In this paper, we prove that Universal Cycles exist for several {\it restricted} classes of words, including non-bijections, ``equitable" words (under suitable restrictions), ranked permutations, and ``passwords".   In each case, proving connectedness of the underlying DeBruijn digraph is the non trivial step.
\end{abstract}

%%%%%%%%%%%%%%%%%%%%%%%%%%%%%%%%%%%%%%%%%%%%%%%%
\section{Introduction}%%%%%%%%%%%%%%%%%%%%%%%%%%%%%%%%%%%%%%%%%%%%%%
The following string has the property that, when wrapped around, it contains all possible words of length three on the binary alphabet \{0,1\}: 11100010.  We call such combinatorial structures universal cycles, or U-cycles, since they list all possible ``values" of a combinatorial object (in this case binary words of length 3) by the mechanism of a sliding window that wraps around the string.  Hurlbert [3] exhibited the following often quoted U-cycle of 3 subsets of \{1,2,...8\}: \[ 1356725\ 6823472\ 3578147\ 8245614\ 5712361\ 2467836\ 7134582\ 4681258, \]
where each block is obtained from the previous one by addition of 5 modulo 8.
Chung, Diaconis and Graham [2], studied U-cycles for permutations, partitions, and $k$-sets of an $n$-set; their work on $k$-subsets of $[n]$ was continued by Hurlbert [3].  More recent work will be featured in a forthcoming issue [8] of {\it Discrete Mathematics} which will contain a selection of papers presented at the Workshop on Generalizations of de Bruijn Cycles and Gray Codes, held at the Banff International Research Station, Banff, Canada, December 4--9, 2004.   Here we will continue the work begun by Bechel et al. [1], who proved that U-cycles of surjections from $\{1,2,\ldots,k\}$ to $\{1,2,\ldots,n\}$ exist iff $n<k$ and who derived a result on ``1-inequitable" binary functions on $\{1,2,\ldots,k\}$, i.e., binary sequences in which the number of zeros and ones differ by exactly one.  Using a variety of proof techniques, we will  offer results on U-cycles for other restricted classes of words, including non-bijections, ``equitable" words, ranked permutations, and ``passwords".

%%%%%%%%%%%%%%%%%%%%%%%%%%%%%%%%%%%%%%%%%%%%%%
\section{Prior Work}%%%%%%%%%%%%%%%%%%%%%%%%%%%%%%%%%%%%%%%%%%%%%%%%
The following result is basic to the theory of U-cycles; see [7] for a proof.
\newtheorem{thm}{Theorem}
\begin{thm} A connected digraph is Eulerian if and only if the in-degree of each vertex is the same as its out-degree.  \end{thm}
\noindent  Theorem 1 is used, for example, to prove the following baseline result on U-cycles, namely de Bruijn's theorem:

\begin{thm} U-cycles of $k$-letter words on an $n$-letter alphabet exist for all $k, n.$  \end{thm} 
\begin{proof}  We create a digraph, $G$, with a vertex set that consists of all $k-1$ letter words on the $n$ letter alphabet.  In other words, vertices have one less letter than the words which we seek to U-cycle, which will appear as edge labels between vertices as follows:  A directed edge is drawn from $v_1$ to $v_2$ if the last $k-2$ letters of $v_1$ are the same as the first $k-2$ letters of $v_2$, and is labelled with the corresponding concatenated $k$-letter word.  For example, the edge from 11234 to 12344 will be labeled 112344.   It is easy to see that the conditions of Theorem 1 are satisfied, and that the Eulerian circuit generates the required U-cycle.  
\end{proof}

\noindent In this paper, we will use a variation of the above proof in all our results.  The key difference is that connectedness of the underlying ``de Bruijn digraph"  is no longer obvious and will need to be proved.  

The following result was proved by Jackson [4].  It shows that a U-cycle of an important class of restricted words exists as well. 

\begin{thm} A U-cycle of 1-1 functions from $\{1,...k \} \rightarrow \{1,...n\}$ exists if and only if $n>k$ ; these are merely permutations of $n$ objects taken $k$ at a time, or, $k$-letter words on $[n]$ in which no letter repeats.  
\end{thm}

When $k=n$, the underlying graph is not connected, and a U-cycle cannot exist. Knuth [5] raised the question of when a U-cycle of one to one functions can be explicitly constructed and the first such effort appears to be, for $k=n-1$, due to Ruskey and Williams [6].   The non-trivial part of the proof of Theorem 3 consists of showing connectedness of the underlying digraph.  The same is true of the next theorem, proved by Bechel et al. [1], who considered words that exhaust the alphabet.

\begin{thm} A U-cycle of onto functions from $\{1,...k\} \rightarrow \{1,...n\}$ exists if and only if $k>n$. 
\end{thm}

In this paper, we prove existence of U-cycles for many more sets of restricted words, such as non-bijections, equitable words, ranked permutations, and ``passwords" (we define these later).  In each case, we have indicated the easiest proof of connectedness that we could find; the proof of Theorem 8 is the most intricate.

%%%%%%%%%%%%%%%%%%%%%%%%%%%%%%%%%%%%%%%%%%%%%
\section{U-cycles of non-bijections} %%%%%%%%%%%%%%%%%%%%%%%%%%%%%%%%%%%%%%%%%%%

Throughout this paper, we will refer to words with entries from an $n$ letter alphabet $\{ 1,2,...n \}$ as being on $[n]$.  
A $k$-letter word on $[n]$, i.e., a function, $f:\{1,...k\} \rightarrow \{1,...n\}$, is said to be \textit{almost onto} if $|[n]- Range(f)|=1$. 
A function $f:[n]\rightarrow[n]$ is said to be a non-bijection if $Range(f)\ne[n]$.

\begin{thm}  A U-cycle of almost onto $n$-letter words on [n] exists for $n\ge3$.\end{thm} 
\begin{proof}  Vertices in our digraph will have labels of $M=n-1$ letter words in which at most one letter appears exactly twice and the other letters are distinct. From now on, we will speak of vertices and the words they represent interchangeably.  For vertices $v$ in which no letter repeats (such vertices will be called $NR$ vertices), $i(v)=o(v)=n-1$, since any letter may be chosen to ``complete the edge label" other than the one which does not appear in $v$; for example with $n=5$, the vertex $1234$ points to 2341, 2342, 2343, and 2344, with the corresponding edge labels being 12341, 12342, 12343, and 12344.  Similarly, it is easy to see that for vertices with exactly  one repeating letter ($OR$ vertices), $i(v)=o(v)=2$, since any one of the two letters which do not appear in $v$ may be chosen to complete the edge label.  In either case, the concatenated word, i.e., the edge label, represents an almost onto word.  Connectedness is easy to establish for $n=3,4$. For $n\ge5$, we shall exhibit a path from an arbitrary $xR$ vertex to another $yR$ vertex, where $x,y$ may equal $N$ or $O$.

\noindent {\it Case 1.}  We prove that it is possible to travel from one $NR$ vertex $A=a_1\ldots a_M$ to another, denoted by $B=b_1\ldots b_M$. There are two steps in the proof.  First we show how we reach a word with the same letters as $B$, and then exhibit an algorithm by which these letters may be put in the ``right order." 

Since $A$ and $B$ are both $NR$, they share all but one letter.  Say $a_j\not\in B; b_i\not\in A$.  Starting with $A$, a word with the same letters as $B$ is reached as follows:
\[a_1\ldots a_j\ldots a_M\rightarrow\ldots\rightarrow a_j\ldots a_M a_1\ldots a_{j-1}\rightarrow a_{j+1}\ldots a_M a_1\ldots a_{j-1}a_{j-1}\rightarrow\]
\[a_{j+2}\ldots a_Ma_1\ldots a_{j-1}a_{j-1}b_i\rightarrow a_{j+3}\ldots a_Ma_1\ldots a_{j-1}a_{j-1}b_ia_{j+1}\rightarrow\ldots\rightarrow\]
\[a_{j-1}a_{j-1}b_ia_{j+1}\ldots a_Ma_1\ldots a_{j-3}\rightarrow a_{j-1}b_ia_{j+1}\ldots a_Ma_1\ldots a_{j-2}.\]
Let the word $a_{j-1}b_ia_{j+1}\ldots a_Ma_1\ldots a_{j-2}$ thus reached be written as $C=c_1\ldots c_M$; $B$ is a permutation of the letters of $C$.  To show that one may travel from $C$ to $B$ it suffices to show that it is possible to go from $C$ to a word created by a single swap, say $D=c_1\ldots c_{i-1}c_jc_{i+1}\ldots c_{j-1}c_ic_{j+1}\ldots c_M$ where $i<j$. Since any permutation is a composition of 2-cycles (``swaps" or ``transpositions") we can thus reach the word $B$. The path from $C$ to $D$ is obtained by (i) cycling till are able to add the first mismatch and thus get a $OR$ word; (ii) adding recently discarded letters until a new $NR$ word is reached; and (iii) repeating the process till the word $D$ is reached.  For example, the path from 12345 to 12543 is as follows:
\[12345\ra23451\ra34512\ra45125\ra51253\ra12534\ra25341\ra\]\[53412\ra34125\ra41254\ra12543.\]

\noindent {\it Case 2.}  To go from a $OR$ word to a $NR$ word, we travel to an intermediate $NR$ word as rapidly as possible and then go from this word to the target $NR$ word as in case 1.

\noindent{\it Case 3.}  If the goal is to describe a path from $A$ ($NR$) to $B$ ($OR$), we identify a $NR$ word $C$ from which $B$ may be reached and then travel from $A$ to $C$ as in Case 1.

\noindent {\it Case 4.} Traveling between two $OR$ words is done by combining Cases 3 and 4.
\end{proof}

\noindent\textbf{Remark} U-cycles of non-bijections on $[n]$ can be shown to exist using the same argument.  We know that for a vertex with no repeats $i(v)=o(v)=n-1$, since we may travel to or from any of the letters that already appear in the vertex.  For a vertex with repeats, $i(v)=o(v)=n$.  The corresponding digraph can easily be shown to be connected:  a path between two $NR$ words is created as in the proof of Theorem 5, and the other 3 cases are similar; for example when traveling from a word with repeats to a $NR$ word, we first eliminate the repeats from the starting word and then travel to the target $NR$ word. 

%%%%%%%%%%%%%%%%%%%%%%%%%%%%%%%%%%%%%%%%%%%%%%%%%%%%%
\section{U-cycles of equitable words}%%%%%%%%%%%%%%%%%%%%%%%%%%%%%%%%%%% 

A word is said to be equitable if for all letters $i,j$, $||i|-|j|| \leq 1$, i.e., if its letters are distributed as evenly as possible.  In this section, we will prove several theorems about U-cycles on equitable vertices.  In Bechel et al. [1], only the binary case was considered.  Here too, the nomenclature differed; words in which the numbers of ones and zeros differed by one were called {\it 1-inequitable}.  We prefer to use the terminology of graph labeling.
%To prove these theorems, we will always assume the worst case.  In this way, we cover all possible problems. 

\begin{thm} A U-cycle of equitable $m$ letter words on $[n]$ where $m \equiv 0 \pmod n$ does not exist. \end{thm} 
\begin{proof} Vertices will have $M=m-1$ letters, and it is clear that $i(v)=o(v)=1\enspace\forall v$;  we must always choose the deficient letter to add on and will end up cycling back to the starting word.  Thus no U-cycle exists. \end{proof}

\bigskip

\noindent {\bf Example:} Here we will take the example of 6 letter words on [3].  Say we have the vertex 11223.  This must go to 12233, to preserve equitability.  Then we must travel to 22331, 23311, 33112, 31122, 11223. We end up where we started.  

\begin{thm} There exists a U-cycle of equitable $m$ letter words on $[n]$, where $m \equiv 1 \pmod n$. \end{thm}
\begin{proof}  Vertices will have $M=m-1$ letters, where $M \equiv 0 \pmod n$.  Let $ r= \frac {M}{n}$.   There are two types of vertices.  In an ``equitable" vertex, there will be $r$ occurrences of each letter.  Such vertices have $i(v)= o(v)= n$, since we can put in any letter to complete the edge label.  ``Inequitable" vertices will have $r$ occurrences of each of $n-2$ letters, with one letter appearing $r-1$ times and the last letter appearing $r+1$ times.  Inequitable  vertices necessarily have $i(v)=o(v)=1$, since we must travel to the word that makes up for the ``deficient" letter. 

We will show connectedness by first exhibiting a path from an arbitrary equitable word $A=a_1 a_2 \ldots a_M$ to another equitable word $B=b_1 b_2\ldots b_M$.  We add letters of the word $B$ until we are legally able to, i.e., until we reach the inequitable vertex $a_s\ldots a_Mb_1\ldots b_{s-1}$. Note that at least one of the $r+1$ occurrences of the ``overrepresented" letter must be among the $a$s, since $B$ is equitable and thus has exactly $r$ letters of each kind; consequently $b_1\ldots b_{s-1}$ contains at most $r$ occurrences of the overrepresented letter.    Let the first occurrence of the overrepresented letter among the $a$s occur at $a_\ell$.  We now add mandatory letters to the inequitable vertex $a_s\ldots a_Mb_1\ldots b_{s-1}$ until we reach the equitable word $a_{\ell+1}\ldots a_Mb_1\ldots b_{s-1}c_1\ldots c_{\ell-s+1}$. This may now be cycled around to get $c_1\ldots c_{\ell-s+1}a_{\ell+1}\ldots a_Mb_1\ldots b_{s-1}$, and, finally, we add $b_s$ to get $c_2\ldots c_{\ell-s+1}a_{\ell+1}\ldots a_Mb_1\ldots b_{s}$.  We repeat the above process until $B$ is reached.  

If we wish to exhibit a path from an inequitable $A$ to an equitable $B$, we first go from $A$ to an equitable $C$ and then from $C$ to $B$ as in the previous paragraph.  The last two cases are handled similarly.
 \end{proof} 

\noindent The next result generalizes Theorem 7; notice, however, that the proof is substantially different.

\begin{thm} A U-cycle of equitable $m$ letter words on $[n]$ exists whenever $m \equiv k \pmod n$, and $k \neq 0$. \end{thm} 
\begin{proof} Assume that $k\ge2$, since the $k=1$ case has already been treated in the previous theorem.  Vertices will have $M=m-1$ letters, where $M \equiv k-1 \pmod n$. Let $r= \frac{m-k}{n}$. Equitable vertices will have $n-k+1$ letters repeated $r$ times, and $k-1$ letters appearing $r+1$ times.  Such vertices will have $i(v)=o(v)=n-k+1$, since we may ``add" any letter of which there are $r$.  Vertices which are inequitable will have $r$ occurrences of $n-k-1$ letters, $r+1$ occurrences of $k$ letters, and $r-1$ occurrences of a single letter.  These vertices will have $i(v)=o(v)=1$, since we must go to the letter that appears $r-1$ times. \\

We will show connectedness by producing a path from an arbitrary equitable word $A=a_1 a_2 \ldots a_M$ to another equitable word $B=b_1 b_2 \ldots b_M$; the other three cases will then follow easily.  Our strategy will be (i) to first exhibit the fact that we can travel from $A$ to a word $B'$ with the same letter frequencies as $B$, and then (ii) to show that we can rearrange the letters of $B'$ as needed; the latter will be done through a series of swaps. 

Letters that appear $r-1$, $r$, and $r+1$ times will be called {\it deficient, normal}, and {\it super}, respectively.  An equitable word has no deficient letters.

\noindent (i) Assume that $A$ has super letters labeled $\a_1,\ldots,\a_{k-1}$ and normal letters $\a_k,\ldots,\a_{n}$.  Similarly let $B$ have super letters $\b_1,\ldots,\b_{k-1}$ and normal letters $\b_k,\ldots,\b_{n}$.  The word frequencies of $B$ are clearly obtained by changing some of the super letters in $A$ to normal letters; each such change forces a normal letter to become a super letter.  A sequence of such ``swaps" permits us to reach the letter frequencies of the target word $B$.  For example, if $A=112233444555666$ and $B=111223334445566$, we need to make `6' normal while converting `1' to super status, and make `5' normal while making `3' super at the same time.  Suppose we wish to ``swap the status" of letters $\a_i;1\le i\le k-1$ and $\a_j;k\le j\le n$ in this fashion.  Suppose furthermore that there are $s\ge1$ super letters that occur before the first occurrence of super letter $\a_i$, and we label these from left to right as $c_1,\ldots,c_s$.  Finally, suppose that there are $t-s$ normal letters, labeled $d_1,\ldots,d_{t-s}$ from left to right, before the first occurrence of $\a_i$. The way in which the $c$s and $d$s are intertwined is irrelevant.  We replace $c_1$ by $\a_j$ so as to make $\a_j$ a super letter; replace each $d_i$ by itself; and restore the super status of $c_1,\ldots,c_s$ by replacing $c_i$ by $c_{i-1}$; $2\le i\le s$ and $\a_i$ by $c_s$.  For example we swap the status of `4' and `1' in the word 512625454331466 by proceeding as follows:
\[512625454331466\rightarrow126254543314661\ra262545433146615\ra\]\[625454331466152\ra254543314661521\ra545433146615212\ra\]\[454331466152126\ra543314661521265.\]If there are no super letters preceding the first $\a_i$, i.e. if $s=0$, we simply replace each normal letter by itself and the first $\a_i$ by $\a_j$.
A sequence of such swaps allows us to arrive at $B'$.

\noindent (ii) We now need to be able to reorder the letters of $B'$ in the order desired, i.e., so as to get $B$. This too will be achieved through a sequence of swaps.  To begin with, however, we show that any equitable word may be ``lag cycled" around if we first add a ``placeholder", which is any normal letter, to the word. Lag cycling is defined to be a process in which one letter is always missing from the cyclic version of the word, and in which there is an extra letter (the placeholder) that is eliminated at the last step.  This initial step makes the cycling legal at all stages, even though an inequitable word may be reached in an intermediate step.  More specifically, if we start with an equitable word $A=a_1\ldots a_M$ containing a normal letter $x$, then the sequence of steps
\[a_1\ldots a_M\ra a_2\ldots a_Mx\ra a_3\ldots a_Mxa_1\ra\ldots\ra xa_1\ldots\a_{M-1}\ra A\]is always legal; for example, if $A=1122333$, the above sequence might be
\[1122333\ra1223332\ra2233321\ra2333211\ra3332112\ra3321122\ra\]\[3211223\ra
2112233\ra1122333.\]The reason that the above process works is evident in hindsight:  The placeholder creates a lag between the deletion of a latter and its re-introduction.  So, e.g., the step $C:=a_2\ldots a_Mx\ra a_3\ldots a_M xa_1=:D$ is always permissible since by design $a_1$ is a normal or deficient letter in the word $C$.  

To be able to swap letters, we are going to need {\it two} placeholders.  First note that since there are $k-1$ super letters and $k\le n-1$, each equitable word {\it must} have at least two normal letters, denoted by $\heartsuit$ and $\spadesuit$, to be used as {placeholders}.  Assume that we need to swap letters $a_i$ and $a_j$, i.e., go from $A=a_1\ldots a_i\ldots a_j\ldots a_M$ to $a_1\ldots a_j\ldots a_i\ldots a_M$.  We start by introducing our first placeholder $\hs$ right away, choosing $\hs=a_i$ if $a_i$ is normal, and $\hs\ne a_j$ if $a_i$ is super:
\[a_1\ldots a_M\ra a_2\ldots a_M\hs,\]
and continue lag cycling as in the previous paragraph until we reach
\[a_{i+1}\ldots a_j\ldots a_M\hs a_1\ldots a_{i-1}.\]
Noting that the above word is equitable, we then introduce the second placeholder $\ss$, {\it and which should be chosen to be $a_j$ if at all possible.} [Note that $a_i$ cannot be deficient at this stage, so we can choose $\ss\ne a_i$. If $a_j$ is a normal letter in the word $A$, we first introduce the other normal letter $\heartsuit$, waiting until this stage to introduce $\ss=a_j$.  Of course if $a_j$ is a super letter then we cannot have $\ss=a_j$ at this stage.]  We thus  get
\[a_{i+2}\ldots a_j\ldots a_M\hs a_1\ldots a_{i-1}\ss,\]
where $\ss$ is as indicated above, and continue lag cycling on to
\[A^*=a_{j+1}\ldots a_M\hs a_1\ldots a_{i-1}\ss a_{i+1}\ldots a_{j-1},\] 
and then, in a critical step, transition to
\[a_{j+2}\ldots a_M\hs a_1\ldots a_{i-1}\ss a_{i+1}\ldots a_{j-1}a_i.\]
Let us check that this last step is valid.  Two things would prevent it from being so.  First, $a_i$ could be super, or, second, $a_j$ could be deficient in the word $A^*$.  
  Consider the first possibility.  Since $\hs\ne\ss$, if neither $\hs$ nor $\ss$ equals $a_i$, then $a_i$ is either normal or deficient in $A^*$.
If $\ss\ne\hs=a_i$, then $a_i$ must have been normal to begin with and still is.  Finally, the possibility $\hs\ne\ss=a_i$ has been ruled out. Can $a_j$ be deficient in $A^*$?  A review of the possibilities shows that this too is impossible.  We thus reintroduce $a_i$ as above into the equitable word $A^*$, and lag cycle to
\[a_M\hs a_1\ldots a_{i-1}\ss a_{i+1}\ldots a_{j-1}a_ia_{j+1}\ldots a_{M-2}.\]
From here, in two steps, and after dropping the first placeholder, we reach
\[A^{**}=a_1\ldots a_{i-1}\ss a_{i+1}\ldots a_{j-1}a_ia_{j+1}\ldots a_{M-1}a_M.\]
If $a_j=\ss$, we are done; the required swap has been achieved.  If $a_j\ne\ss$, $a_j$ is a normal letter in the equitable word $A^{**}$ and we add another available placeholder $\cs$ to get 
\[a_2\ldots a_{i-1}\ss a_{i+1}\ldots a_{j-1}a_ia_{j+1}\ldots a_{M-1}a_M\cs,\]
and lag cycle until we reach
\[A^{***}=a_{i+1}\ldots a_{j-1}a_ia_{j+1}\ldots a_M\cs a_1\ldots a_{i-1}.\]
Next, we reintroduce the (now) normal letter $a_j$ into the equitable word $A^{***}$ to yield
\[a_{i+2}\ldots a_{j-1}a_ia_{j+1}\ldots a_M\cs a_1\ldots a_{i-1}a_j,\]
and lag cycle until the target word 
\[a_1\ldots a_{i-1}a_j a_{i+1}\ldots a_{j-1}a_ia_{j+1}\ldots a_{M-1}a_M\]
is reached.

\noindent 

\end{proof}

\noindent {\bf Remark} Suppose we define an {\it s-inequitable} word as one in which letter frequencies differ by at most $s$, i.e., $\vert\vert i\vert-\vert j\vert\vert\le s$ for all $i,j$; $s\ge 1$.  (Using this nomenclature, equitable words could be termed {\it 1-inequitable}, as was done in in Bechel et al. [1]). It is then easy to use Theorem 8 to show that a U-cycle of $s$ inequitable $m$ letter words exists as long as $m$ is not a multiple of $n$.  Here is a sketch of the proof:  Assume for simplicity that $s=2$.  Vertices are words of length $M=m-1$ and may be of three kinds -- they may be 1-inequitable (or, equitable), 2-inequitable, or 3-inequitable. For example, if $n=5$ and $m=19$, then vertices may have letter frequencies $(4,4,4,3,3)$, or $(4,4,4,4,2)$, or $(5,4,4,3,2)$. In general there are as many cases as there are partitions of the integer $m-1$ into $n$ parts with the difference between the maximum part size and the minimum part size being at most $s+1$, and with the minimum size part size having multiplicity one when the above difference equals $s+1$. In our case $s+1=3$.  For these three types of vertices, $i(v)=o(v)=5, 1, 1$ respectively.  To establish connectedness, we first travel from the ``start" word to a 1-inequitable one, note than we can backtrack from the target word to another 1-inequitable word, and, finally, go from the first 1-inequitable word thus created to the second as in the proof of Theorem 8. 
%%%%%%%%%%%%%%%%%%%%%%%%%%%%%%%%%%%%%%%%%

\section{U-cycles of Ranked Permutations and Passwords}%%%%%%%%%%%%%%%%%%%%%%%%%

Chung, Diaconis and Graham [2] suggest investigating U-cycles on tied permutations as a future direction of research.  Here we offer one way of defining these, motivated by rankings, seedings, etc. in sports and other events. 

\medskip

\noindent {\bf Definition} We say that a word satisfies a ranking if the word contains a 1, and if there exists $r\ge1$ of some letter $a$, the next letter must be $a+r$. 

More informally, these words must follow ordinary rankings in a tournament.  For example, the ranking 113 is allowed, but not 112, since second place is already taken in the tie for first.

\medskip

\noindent {\bf Example} Here are all the rankings allowed on $[3]$: 
\[111,122,221,212,113,311,131,123,312,231,321,132,213.\] We see that a U-cycle of these words exists as follows: 1113212213123.  Of the associated vertices, 22, 23, and 32 have $i(v)=o(v)=1$ and 12, 21, 31, 13, 11 have $i(v)=o(v)=2$.  An enumeration of rankings on $[n]$ for $1\le n\le 18$ (the so-called ordered Bell numbers) may be found as sequence A000670 in Neil Sloane's website of integer sequences 
\[{\tt www.research.att.com/}\sim{\tt njas/sequences/}\]

\begin{thm} A U-cycle of ranked permutations on $m$-letter words exists for each $m$. \end{thm}
\begin{proof}  Vertices in the graph consist of $m-1$ letter words that are consistent with a ranking, i.e., those that can be extended to a ranking.  For example, the word 1124 is not consistent with a ranking, but 1135 is.  Vertices will either contain a 1 or not.  Since all rankings must contain a 1, vertices which do not contain a 1 will have $i(v)=o(v)=1$.  Consider vertices which do contain a 1.  Now each vertex will be missing one letter from some ranking.  We now ask: how many rankings is such a vertex consistent with?  We claim the answer is 2.  Writing a vertex with its letters in non decreasing order, we see that there is a single ``gap" in a ranking, and that gap may always be filled in in two ways -- by the beginning letter of the gap or the next possible letter.  (For example, say we had 5 letter words, and we had the vertex 2125.  This is really the ranking $122x5$ where the $x$ may be filled by another 2 or by a 4, the next letter in a logical ranking. The vertex 1114557899 has a gap at the end which may be filled with a 9 or an 11.) Thus $i(v)=o(v)=2$.

Connectedness:  From any vertex we can travel in the direction of less repeats, until we get to a vertex without any ties.   Likewise, we can backtrack from any vertex to one without repeats. It remains to be shown that we can travel from any one vertex $A$ without repeats to another, labeled $B$, that also has no repeats.  Our strategy, distinct from that adopted in previous proofs, will be to show that we may legally travel from $A$ to $O=111\ldots 111$ and then from $O$ to $B$.  Traveling from $A$ to $O$:  Since $A$ has no repeats, we first indicate how to replace the `2', it it exists, by a `1', transitioning in the process to a word $A_1$ with two 1s.  We add letters to $A$ until the 2 is eliminated.  We may now add a 1, and may need to add another 1 in order for the word to have two 1s, as desired.  As an illustration, we implement the replacement of 2 by 1 as follows: 
\[A:=532147\ra321476\ra214765\ra147653\ra476531\ra765311:=A_1.\]
We next add letters to $A_1$ until the `3', if it exists, is eliminated - and then add an additional 1 in at most three steps, thus getting a word $A_2$.  Continuing in this fashion, we reach $O$; in the above example, our procedure would unfold as follows:
\[A_1=765311\ra653114\ra531147\ra311476\ra114765\ra147651\ra\]\[476511\ra765111=:A_2,\]
and be completed thus:
\[765111\ra651111\ra511117\ra111176\ra\ldots\ra761111\ra611111\ra\]
\[111117\ra\ldots\ra711111\ra111111=O.\]
Traveling from $O$ to $B$ is essentially a reversal of the above process.  Given a word $B$, define $B^*$ to be the word $B$ read backwards.  Let $P^*$ be a path from $B^*$ to $O$ created as in the previous paragraph.  Then the path $P$ defined to be $P^*$ read backwards, and with entries read backwards too, legally takes us from $O$ to $B$.  For example, to go from $O$ to 741235 (the starting word in the previous paragraph read backwards), we use the steps
\[O\ra111117\ra\ldots\ra111567\ra\dots\ra113567\ra\ldots\ra\]\[567412\ra674123\ra741235.\]\end{proof} 

\medskip

\noindent{\bf Definition} We say that an $m$-letter word on $[n]$ is a {\it password} if there are $q<n$ distinct classes of symbols in $n$, and each word must contain at least one element of each class.  (Note that the classes need not form a partition of $[n]$.)

More informally, we can think of this as being like a security-conscious Internet password that must contain one lower case letter, one number, one symbol, and so on.  This general definition of passwords includes more familiar elementary textbook objects such as words with at least one vowel, etc. Note moreover that when $m>n=q$, passwords are onto functions; and when $q=1$ we get the set up of DeBruijn's theorem.  Finally we observe that U-cycles of passwords cannot possibly exist for all values of the parameters -- e.g., if $n=m=q$, when passwords are $n$-permutations of $[n]$.  Our next result can almost certainly be improved, and we invite the reader to do so.

\begin{thm}A U-cycle exists for all $m$-letter passwords on $[n]$, with $q$ distinct classes of symbols, provided that $m\ge 2q$.  \end{thm} 
\begin{proof} Vertices will have $M=m-1$ letters.  There are two types of vertices.  If a vertex is missing exactly one of the types of symbols, then $i(v)=o(v)=|l_v|$, where $|l_v|$ is the number of the type of missing symbol, $l_v$, that the vertex is missing.  If a vertex has all the different types of symbols, then clearly $i(v)=o(v)=n$.

We will show connectedness by going from an arbitrary word 

\noindent $A=a_1 a_2 \ldots a_M$, to a word $B=b_1 b_2 \ldots b_M$.  As with previous proofs, the main case (among four altogether) is the one where both $A$ and $B$ have representatives from all classes of symbols.  We first go from $A$ to $A'$, where $A'$ has one of each type of special symbol as its last $q$ letters. Next, we travel from $A'$ to $B'$, where $B'$ has as its last $q$ letters one of each of the special classes of symbols, but in the same order as they appear in the word $B$.  The word $B$ may now be built without hindrance, one letter at a time.  \end{proof}

\section{Acknowledgment}  The research of both authors was supported by NSF REU Grant 0552730, and conducted at East Tennessee State University in the Summer of 2008, when Leitner was a student at California State University, Chico

%%%%%%%%%%%%%%%%%%%%%%%%%%%%%%%%%%%%%%%%%%%%%%
\section{References}

[1]  A.~Bechel, B.~LaBounty-Lay and A.~Godbole (2008), ``Universal cycles of discrete functions,"  {\it Congressus Numerantium} {\bf 189}, 121--128. 

\noindent [2]  F.~Chung, P.~Diaconis, and R.~Graham (1992), ``Universal Cycles for combinatorial structures," \textit{Discrete Math.} \textbf{110}, 43-59. 

\noindent [3]  G.~Hurlbert (1994), ``On universal cycles for $k$-subsets of an $n$-element set," \textit{ SIAM J. Discrete Math}. \textbf{7}, 598-604.

\noindent [4]  B.~Jackson (1993), ``Universal cycles of $k$-subsets and $k$-permutations," \textit{Discrete Math}. \textbf{117}, 114-150.  

\noindent [5] D.~Knuth (2005), {\it The Art of Computer Programming, Volume 4, Fascicle 2}, Pearson, NJ.

\noindent [6] F.~Ruskey and A.~Williams (2008), ``An explicit universal cycle for the $n-1$-permutations of an $n$-set," Talk at Napier Workshop.

\noindent [7]  D.~West  (1996), \textit{Introduction to Graph Theory}, Prentice Hall, New Jersey. 

\noindent [8] Forthcoming Special Issue of {\it Discrete Mathematics}, Glenn Hurlbert and Brett Stevens, eds.

% Set the ending of a LaTeX document
\end{document}